\documentclass{article}%
\usepackage{amsmath}
\usepackage{amsfonts}
\usepackage{amssymb}
\usepackage{color}
\usepackage{graphicx}%
\providecommand{\U}[1]{\protect\rule{.1in}{.1in}}
\newtheorem{theorem}{Theorem}

\newtheorem{lemma}[theorem]{Lemma}

\newtheorem{proposition}[theorem]{Proposition}

\newenvironment{proof}[1][Proof]{\noindent\textbf{#1.} }{\ \rule{0.5em}{0.5em}}

\begin{document}

\title{Functional Estimates for Derivatives \\of the Modified Bessel Function $K_{0}$ and related Exponential Functions}
\author{Silvia Falletta\thanks{Dip. Scienze Matematiche \textquotedblleft G.L.
Lagrange\textquotedblright, Politecnico di Torino, C.so Duca degli Abruzzi 24,
10129 Torino, Italy, e-mail: \texttt{silvia.falletta@polito.it}}
\and Stefan A. Sauter\thanks{Institut f\"{u}r Mathematik, Universit\"{a}t
Z\"{u}rich, Winterthurerstrasse 190, CH-8057 Z\"{u}rich, Switzerland, e-mail:
\texttt{stas@math.uzh.ch}}}
\maketitle

\begin{abstract}
Let $K_{0}$ denote the modified Bessel function of second kind and zeroth
order. In this paper we will studying the function $\tilde{\omega}_{n}\left(
x\right)  :=\frac{\left(  -x\right)  ^{n}K_{0}^{\left(  n\right)  }\left(
x\right)  }{n!}$ for positive argument. The function $\tilde{\omega}_{n}$
plays an important role for the formulation of the wave equation in two
spatial dimensions as a retarded potential integral equation. We will prove
that the growth of the derivatives $\tilde{\omega}_{n}^{\left(  m\right)  }$
with respect to $n$ can be bounded by $O\left(  \left(  n+1\right)
^{m/2}\right)  $ while for small and large arguments $x$ the growth even
becomes independent of $n$.

These estimates are based on an integral representation of $K_{0}$ which
involves the function $g_{n}\left(  t\right)  =\frac{t^{n}}{n!}\exp\left(
-t\right)  $ and their derivatives. The estimates then rely on a subtle
analysis of $g_{n}$ and its derivatives which we will also present in this paper.

\end{abstract}

\noindent\textbf{Keywords:} Bessel Function $K_{0}$ of second kind and order
zero, exponential function, functional error estimates.

\noindent\textbf{Mathematics Subject Classification (2000):} 26A12, 39B62,
33C10, 65R20

\section{Introduction}

The study of Bessel functions is a classical field in mathematics and a vast
literature is devoted to its analysis. Classical references include
\cite[first edition: 1922.]{Watson}, \cite{Abramowitz}, \cite[Chap.
9]{Temme_bessel}, \cite{Olver_bessel}, where the main focus is on the detailed
study of the asymptotic behavior in various parameter regimes, functional
relations, and recurrences.

Sharp estimates of quotients of modified Bessel functions are proved, e.g., in
\cite{Baricz_bessel}, \cite{Ifantis_bessel}, \cite{Joshi_bessel}, while
Tur\'{a}n-type estimates can be found in \cite{Baricz_turan}. As is common in
applied mathematics new questions arise from new applications; in this case,
the rapidly increasing interest in \textit{convolution quadrature} methods for
solving the wave equation in exterior domains (cf.
\cite{Lub1,Lub2,LubOstermann,Lub3,lopezsauter_paper1,BanjaiGruhne}) leads to
the question of functional-type inequalities for derivatives of the Bessel
functions which are explicit with respect to the order of the derivatives.
This has been invested for the three-dimensional case in \cite{ks3}, where the
kernel function is the exponential function.

The modified Bessel function $K_{0}$ plays an important role when studying the
wave equation in two spatial dimensions as retarded potential integral
equations. Their fast numerical solution with panel-clustering and multipole
methods requires the approximation of these functions with polynomials. In
this paper, we will prove new and sharp functional estimates for the function
$\tilde{\omega}_{n}\left(  x\right)  :=\frac{\left(  -x\right)  ^{n}%
K_{0}^{\left(  n\right)  }\left(  x\right)  }{n!}$ and its derivatives on the
positive real axes. It turns out that the growth of the $m$-th derivative
$\tilde{\omega}_{n}^{\left(  m\right)  }$ can be bounded, in general, by
$O\left(  \left(  n+1\right)  ^{m/2}\right)  $. In addition, we will prove
that for small and large argument $x$ the derivatives can be bounded
independent of $n$. This is at first glance surprising since a straightforward
use of Cauchy's integral theorem leads to an \textit{exponential }growth with
respect to $n$. Our proof is based on an integral representation of $K_{0}$,
where the function $g_{n}\left(  t\right)  =\frac{t^{n}}{n!}\exp\left(
-t\right)  $ is involved. The estimates then rely on a subtle analysis of
$g_{n}$ and its derivatives which we will also present in this paper.

These results imply that fast panel-clustering techniques can efficiently be
used to solve retarded potential integral equations in two spatial dimensions.

In \cite{BanjaiGruhne}, cutoff techniques have been introduced to accelerate
the computation of two dimensional retarded potential equations. For this
purpose, the behavior of the function $\omega_{n}$ for large arguments has
been investigated which leads to a speed-up for large times compared to
conventional methods.

\section{Functional Estimates for Derivatives of the Modified Bessel Function}

The kernel function for acoustic retarded potential integral equations --
discretized by convolution quadrature with the $\operatorname*{BDF1}$ method
-- is related to the function $\tilde{\omega}_{n}\left(  x\right)
:=\frac{\left(  -x\right)  ^{n}K_{0}^{\left(  n\right)  }\left(  x\right)
}{n!}$ by%
\[
\omega_{n}\left(  d\right)  =\frac{1}{2\pi}\tilde{\omega}_{n}\left(  \frac
{d}{\Delta t}\right)  .
\]

The estimate of its local approximability by polynomials requires the
investigation of the $m$-th derivative of $\tilde{\omega}_{n}$.

\begin{theorem}
\label{Lemdomegan}

\begin{enumerate}
\item General estimate. For all $n,m\in\mathbb{N}_{0}$, and $x>0$, the
estimate%
\begin{equation}
\left\vert \tilde{\omega}_{n}^{\left(  m\right)  }\left(  x\right)
\right\vert \leq m!\left(  \frac{\gamma\sqrt{n+1}}{x}\right)  ^{m}\left(
\frac{\delta_{m,0}}{\sqrt{2x+1}}+\frac{\gamma}{\sqrt{n+1}}\left(
1+\delta_{m+n,0}\log\frac{x+1}{x}\right)  \right)  \label{estomegatildenm1}%
\end{equation}
holds for some $\gamma\geq1$ (independent of $m$, $n$, and $x$); $\delta
_{m,n}$ denotes the Kronecker delta.

\item Refined estimates for small and large arguments.

\begin{enumerate}
\item Small argument. For given $C>1$ independent of $m,n,x$, there exists
some constant $\gamma>1$ independent of $m,n,x$ such that for all $n\geq0$ and
$0\leq m\leq\sqrt{n}\left(  \frac{C-1}{C}\right)  $ with the further
restriction on $m$:%
\begin{equation}
2\left\lfloor \frac{m\log\left(  n+1\right)  }{4\log2}\right\rfloor
-2\leq\frac{n-1}{2} \label{addcondm}%
\end{equation}
and all
\begin{equation}
0<x\leq\min\left\{  \frac{n}{2C},\frac{n+1}{4\gamma}\right\}
\label{condxsmallargmts}%
\end{equation}
it holds%
\begin{equation}
\left\vert \tilde{\omega}_{n}^{\left(  m\right)  }\left(  x\right)
\right\vert \leq m!\left(  \frac{\gamma}{x}\right)  ^{m}\left(  \frac
{\delta_{m,0}}{\sqrt{2x+1}}+\frac{\gamma}{\sqrt{n+1}}\left(  1+\delta
_{n+m,0}\log\frac{x+1}{x}\right)  \right)  . \label{estomegatildenm1sa}%
\end{equation}

\item Large argument. For all $n\geq0$ and $m\geq1+2\log\left(  n+1\right)  $
it holds%
\begin{equation}
\left\vert \tilde{\omega}_{n}^{\left(  m\right)  }\left(  x\right)
\right\vert \leq m!\left(  \frac{\gamma}{x}\right)  ^{m}\qquad\forall
x>\left\{
\begin{array}
[c]{ll}%
0 & n=0,1,\\
n+m\left(  \sqrt{n}+2\right)  & n\geq2,
\end{array}
\right.  \label{estomegatildenm1la}%
\end{equation}
for some constant $\gamma\geq1$ (independent of $m$, $n$, and $x$).
\end{enumerate}

\item Exponential decay. For $m=0$ and $x\geq\max\left\{  1,n+\sqrt
{n}\right\}  $, the function $\tilde{\omega}_{n}$ is decaying exponentially%
\begin{equation}
\left\vert \tilde{\omega}_{n}\left(  x\right)  \right\vert \leq3\frac
{\exp\left(  \sqrt{n}-\frac{x}{1+\sqrt{n}}\right)  }{\sqrt{n+1}}.
\label{exponentialdecayomegatilde}%
\end{equation}

\end{enumerate}
\end{theorem}%

%TCIMACRO{\TeXButton{Proof}{\proof}}%
%BeginExpansion
\proof
%EndExpansion
We have%
\[
\tilde{\omega}_{n}\left(  x\right)  :=\frac{\left(  -x\right)  ^{n}%
K_{0}^{\left(  n\right)  }\left(  x\right)  }{n!}\overset{\text{\cite[9.6.23]%
{Abramowitz}}}{=}\frac{x^{n}}{n!}\int_{1}^{\infty}\frac{t^{n}\operatorname*{e}%
\nolimits^{-xt}}{\sqrt{t^{2}-1}}dt=\frac{1}{n!}\int_{0}^{\infty}\frac{\left(
x\left(  s+1\right)  \right)  ^{n}\operatorname*{e}\nolimits^{-x\left(
s+1\right)  }}{\sqrt{s\left(  s+2\right)  }}ds.
\]
Let $g_{n}\left(  t\right)  :=\frac{t^{n}\operatorname*{e}\nolimits^{-t}}{n!}%
$. The $m$-th derivative of $\tilde{\omega}_{n}$ is given by
\begin{align*}
\tilde{\omega}_{n}^{\left(  m\right)  }\left(  x\right)   &  =\frac{1}{n!}%
\int_{0}^{\infty}\frac{\left(  \frac{d}{dx}\right)  ^{m}\left(  \left(
x\left(  s+1\right)  \right)  ^{n}\operatorname*{e}\nolimits^{-x\left(
s+1\right)  }\right)  }{\sqrt{s\left(  s+2\right)  }}ds\\
&  =\int_{0}^{\infty}\frac{\left(  s+1\right)  ^{m}g_{n}^{\left(  m\right)
}\left(  x\left(  s+1\right)  \right)  }{\sqrt{s\left(  s+2\right)  }}ds\\
&  =\frac{1}{x^{m}}\int_{x}^{\infty}\frac{t^{m}g_{n}^{\left(  m\right)
}\left(  t\right)  }{\sqrt{t^{2}-x^{2}}}dt.
\end{align*}
We split this integral into $\tilde{\omega}_{n}^{\left(  m\right)  }%
=P_{n,L}^{\left(  m\right)  }+Q_{n,L}^{\left(  m\right)  }$ for some $L\geq x$
with%
\[
P_{n,L}^{\left(  m\right)  }\left(  x\right)  :=\frac{1}{x^{m}}\int_{x}%
^{L}\frac{t^{m}g_{n}^{\left(  m\right)  }\left(  t\right)  }{\sqrt{t^{2}%
-x^{2}}}dt\quad\text{and\quad}Q_{n,L}^{\left(  m\right)  }\left(  x\right)
:=\frac{1}{x^{m}}\int_{L}^{\infty}\frac{t^{m}g_{n}^{\left(  m\right)  }\left(
t\right)  }{\sqrt{t^{2}-x^{2}}}dt
\]
and introduce the quantity
\[
M_{a,b}^{\ell,m,n}:=\sup_{a\leq t\leq b}\left\vert t^{m+\ell}g_{n}^{\left(
m\right)  }\left(  t\right)  \right\vert .
\]
\textbf{Estimate of} $P_{n,L}^{\left(  m\right)  }$\textbf{.}

For general $m,n\geq0$ and $0<x\leq L$ we get%

\begin{align}
\left\vert P_{n,L}^{\left(  m\right)  }\left(  x\right)  \right\vert  &
\leq\frac{M_{x,L}^{0,m,n}}{x^{m}}\int_{x}^{L}\frac{1}{\sqrt{t^{2}-x^{2}}%
}dt\nonumber\\
&  =\frac{M_{x,L}^{0,m,n}}{x^{m}}\left(  \log\frac{L}{x}+\log\left(
1+\frac{\sqrt{L^{2}-x^{2}}}{L}\right)  \right)  \leq\frac{M_{x,L}^{0,m,n}%
}{x^{m}}\log\frac{2L}{x}. \label{estpnm0}%
\end{align}
For our final estimate, we will choose $L\in\left\{  x+1,2x\right\}  $. For
$L=2x$, we get%
\begin{equation}
\left\vert P_{n,2x}^{\left(  m\right)  }\left(  x\right)  \right\vert
\leq\frac{M_{x,2x}^{0,m,n}}{x^{m}}\log4. \label{PnmL2x}%
\end{equation}
For $L=x+1$, we have to refine the estimate to see that the logarithmic
behavior of $P_{n,L}^{\left(  m\right)  }\left(  x\right)  $ as $x\rightarrow
0$ only appears for $n=m=0$. We obtain

\begin{itemize}
\item for $x\geq1$ directly from (\ref{estpnm0}):%
\begin{equation}
\left\vert P_{n,x+1}^{\left(  m\right)  }\left(  x\right)  \right\vert
\leq\frac{M_{x,x+1}^{0,m,n}}{x^{m}}\log4; \label{PnmLxp1la}%
\end{equation}

\item for $0\leq x\leq1$

\begin{itemize}
\item[i.] for $n=m=0$%
\begin{equation}
\left\vert P_{0,x+1}^{\left(  0\right)  }\left(  x\right)  \right\vert
\leq\log\frac{2\left(  x+1\right)  }{x}; \label{PnmLxp1sa1}%
\end{equation}

\item[ii.] for $n=0\wedge m\geq1$%
\begin{align}
\left\vert P_{0,x+1}^{\left(  m\right)  }\left(  x\right)  \right\vert  &
\leq\dfrac{\left(  m-1\right)  !}{x^{m}}M_{x,x+1}^{0,0,m-1}%
%TCIMACRO{\dint _{x}^{x+1}}%
%BeginExpansion
{\displaystyle\int_{x}^{x+1}}
%EndExpansion
\dfrac{t}{\sqrt{t^{2}-x^{2}}}dt\nonumber\\
&  =\sqrt{3}\dfrac{\left(  m-1\right)  !}{x^{m}}M_{x,x+1}^{0,0,m-1};
\label{PnmLxp1sa2}%
\end{align}

\item[iii.] for $n\geq1\wedge m=0$%
\begin{equation}
\left\vert P_{n,x+1}^{\left(  0\right)  }\left(  x\right)  \right\vert
\leq\int_{x}^{x+1}\frac{t^{n}\operatorname*{e}^{-t}}{n!\sqrt{t^{2}-x^{2}}%
}dt\leq\sqrt{3}\frac{M_{x,x+1}^{0,0,n-1}}{n}; \label{PnmLxp1sa3}%
\end{equation}

\item[iv.] for $n\geq1\wedge m\geq1$

By employing the recursion $g_{n}^{\left(  m\right)  }=g_{n-1}^{\left(
m-1\right)  }-g_{n}^{\left(  m-1\right)  }$ we get%
\begin{align}
\left\vert P_{n,x+1}^{\left(  m\right)  }\left(  x\right)  \right\vert  &
\leq\frac{M_{x,x+1}^{0,m-1,n-1}+M_{x,x+1}^{0,m-1,n}}{x^{m}}\int_{x}^{x+1}%
\frac{t}{\sqrt{t^{2}-x^{2}}}dt\nonumber\\
&  =\sqrt{3}\frac{M_{x,x+1}^{0,m-1,n-1}+M_{x,x+1}^{0,m-1,n}}{x^{m}}.
\label{PnmLxp1sa4}%
\end{align}

\end{itemize}
\end{itemize}

For $m=0$ and $x\geq n+\sqrt{n}$, the estimate can be refined by using
(\ref{estgnasymptotic}):%
\begin{align}
\left\vert P_{n,L}^{\left(  0\right)  }\left(  x\right)  \right\vert  &
\leq\int_{x}^{L}\frac{g_{n}\left(  t\right)  }{\sqrt{t^{2}-x^{2}}}%
dt\leq\left(  \max_{x\leq t\leq L}g_{n}\left(  t\right)  \right)  \log
\frac{2L}{x}\nonumber\\
&  \leq\frac{\log\frac{2L}{x}}{\sqrt{n+1}}\exp\left(  \sqrt{n}-\frac
{x}{1+\sqrt{n}}\right)  . \label{estpn0}%
\end{align}

\textbf{Estimate of }$Q_{n,L}^{\left(  m\right)  }\left(  x\right)
$\textbf{.}

For $m=0$ it is easy to see that
\begin{equation}
Q_{n,L}^{\left(  0\right)  }\left(  x\right)  =\int_{L}^{\infty}\frac
{g_{n}\left(  t\right)  }{\sqrt{t^{2}-x^{2}}}dt\leq\frac{1}{\sqrt{L^{2}-x^{2}%
}}\int_{0}^{\infty}g_{n}\left(  t\right)  dt=\frac{1}{\sqrt{L^{2}-x^{2}}}.
\label{estqn00}%
\end{equation}
In our application, we will employ the function $Q_{n,L}^{\left(  0\right)  }$
only for $L:=x+1$ and obtain from (\ref{estqn00})%
\begin{equation}
Q_{n,x+1}^{\left(  0\right)  }\left(  x\right)  \leq\frac{1}{\sqrt{2x+1}}.
\label{est2xp1}%
\end{equation}
Again, for $x\geq n+\sqrt{n}$, this estimate can be refined by using
(\ref{estgnasymptotic}):%

\begin{equation}
Q_{n,x+1}^{\left(  0\right)  }\left(  x\right)  \leq\frac{1}{\sqrt{2x+1}}%
\frac{1}{\sqrt{n+1}}\int_{L}^{\infty}\exp\left(  \sqrt{n}-\frac{t}{1+\sqrt{n}%
}\right)  dt\leq\sqrt{2}\frac{\exp\left(  \sqrt{n}-\frac{x}{1+\sqrt{n}%
}\right)  }{\sqrt{2x+1}}. \label{estqn0}%
\end{equation}

Next we consider $m\geq1$ and introduce a splitting $Q_{n,L}^{\left(
m\right)  }=R_{n,L,\mu}^{\left(  m\right)  }+S_{n,L,\mu}^{\left(  m\right)  }$
which we will explain next. First, note that the Taylor expansion of $\left(
t^{2}-x^{2}\right)  ^{-1/2}$ around $x=0$ can be written in the form (cf.
\cite[(5.24.31)]{Hansen75})%
\begin{equation}
\frac{1}{\sqrt{t^{2}-x^{2}}}\approx\frac{T_{\mu}\left(  x/t\right)  }{t}%
\quad\text{with\quad}T_{\mu}\left(  z\right)  :=\sum_{\ell=0}^{\mu-1}%
\tbinom{2\ell}{\ell}\left(  \frac{z}{2}\right)  ^{2\ell}\qquad\mu\in
\mathbb{N}_{\geq1}. \label{defTmue}%
\end{equation}
We set%
\begin{align*}
R_{n,L,\mu}^{\left(  m\right)  }\left(  x\right)   &  :=\frac{1}{x^{m}}%
\int_{L}^{\infty}t^{m}g_{n}^{\left(  m\right)  }\left(  t\right)  \left(
\frac{1}{\sqrt{t^{2}-x^{2}}}-\frac{T_{\mu}\left(  x/t\right)  }{t}\right)
dt,\\
S_{n,L,\mu}^{\left(  m\right)  }\left(  x\right)   &  :=\frac{1}{x^{m}}%
\int_{L}^{\infty}T_{\mu}\left(  x/t\right)  t^{m-1}g_{n}^{\left(  m\right)
}\left(  t\right)  dt.
\end{align*}

For the estimate of $R_{n,L,\mu}^{\left(  m\right)  }$, we employ Proposition
\ref{PropTaylor} to obtain%

\begin{equation}
R_{n,L,\mu}^{\left(  m\right)  }\left(  x\right)  \leq\frac{M_{L,\infty
}^{0,m,n}}{x^{m}}\int_{L}^{\infty}\frac{1}{\sqrt{t^{2}-x^{2}}}\left(  \frac
{x}{t}\right)  ^{2\mu}dt\leq\frac{M_{L,\infty}^{0,m,n}}{x^{m}}\left(  \frac
{x}{L}\right)  ^{2\mu}. \label{EstRnm}%
\end{equation}

\textbf{Estimate of }$S_{n,L,\mu}^{\left(  m\right)  }$ \textbf{for} $m\geq1$
\textbf{and} $n=0.$

By using
\[
\left\vert T_{\mu}\left(  z\right)  \right\vert \leq\sum_{\ell=0}^{\mu
-1}\tbinom{2\ell}{\ell}\left(  \frac{1}{2}\right)  ^{2\ell}=2^{1-2\mu}%
\mu\binom{2\mu}{\mu}=:c_{\mu}\qquad\forall z\in\left[  0,1\right]
\]
the function $S_{n,L,\mu}^{\left(  m\right)  }$ can be estimated for $n=0$ by%
\begin{equation}
\left\vert S_{0,L,\mu}^{\left(  m\right)  }\left(  x\right)  \right\vert
\leq\frac{c_{\mu}}{x^{m}}\int_{0}^{\infty}t^{m-1}\operatorname*{e}%
\nolimits^{-t}dt=c_{\mu}\frac{\left(  m-1\right)  !}{x^{m}}. \label{estsnmm0}%
\end{equation}

\textbf{Estimate of }$S_{n,L,\mu}^{\left(  m\right)  }$ \textbf{for} $m\geq1$
\textbf{and} $n\geq1$

First, observe that for $0\leq2r\leq n-1$ it holds (cf. Proposition
\ref{Propgnderv})%
\[
\int_{L}^{\infty}t^{m-1-2r}g_{n}^{\left(  m\right)  }\left(  t\right)
dt=-\frac{\left(  n-1-2r\right)  !}{n!}\left.  \left(  y^{m}g_{n-1-2r}%
^{\left(  m-1-2r\right)  }\left(  y\right)  \right)  ^{\left(  2r\right)
}\right\vert _{y=L}.
\]
Note that%
\begin{align}
\left\vert \left(  y^{m}g_{n-1-2r}^{\left(  m-1-2r\right)  }\left(  y\right)
\right)  ^{\left(  2r\right)  }\right\vert _{y=L}  &  =\left\vert \sum
_{\ell=0}^{\min\left(  m,2r\right)  }\tbinom{2r}{\ell}\frac{m!}{\left(
m-\ell\right)  !}L^{m-\ell}g_{n-1-2r}^{\left(  m-1-\ell\right)  }\left(
L\right)  \right\vert \nonumber\\
&  \leq\sum_{\ell=0}^{\min\left(  m,2r\right)  }\tbinom{2r}{\ell}\frac
{m!}{\left(  m-\ell\right)  !}M_{L,L}^{1,m-1-\ell,n-1-2r}=:K_{n,L}^{m,r}.
\label{DefKL0}%
\end{align}
Thus,%
\begin{align}
\left\vert S_{n,L,\mu}^{\left(  m\right)  }\left(  x\right)  \right\vert  &
\leq\frac{1}{x^{m}}\sum_{r=0}^{\mu-1}\left(  \frac{x}{2}\right)  ^{2r}%
\tbinom{2r}{r}\left\vert \int_{L}^{\infty}t^{m-1-2r}g_{n}^{\left(  m\right)
}\left(  t\right)  dt\right\vert \nonumber\\
&  \leq\frac{1}{x^{m}}\sum_{r=0}^{\mu-1}\left(  \frac{x}{2}\right)
^{2r}\tbinom{2r}{r}\frac{\left(  n-1-2r\right)  !}{n!}K_{n,L}^{m,r}.
\label{estsnm0}%
\end{align}

The derived estimates imply the assertion by arguing as follows.

\textbf{Case 1: Exponential decay.}

Estimate (\ref{exponentialdecayomegatilde}) follows from (\ref{estpn0}) and
(\ref{estqn0}) (with $x\geq1$ and $L:=x+1$).\medskip

\textbf{Case 2: General estimate and estimate for large argument.}

For $m=0$, estimate (\ref{PnmLxp1sa1}), (\ref{PnmLxp1la}), (\ref{PnmLxp1sa3}),
(\ref{est2xp1}), (\ref{estgntilde}) imply with $L=x+1$%
\begin{equation}
\left\vert \tilde{\omega}_{n}\left(  x\right)  \right\vert \leq\frac{1}%
{\sqrt{2x+1}}+\frac{1}{\sqrt{n+1}}\left\{
\begin{array}
[c]{ll}%
\log\frac{4\left(  x+1\right)  }{x} & n=0,\\
\gamma & n\geq1,
\end{array}
\right.  \label{gnfin1}%
\end{equation}
for some $\gamma>1$ so that (\ref{estomegatildenm1}) and
(\ref{estomegatildenm1la}) follow for $m=0$.

For $m\geq1$, $n=0$, the estimates (\ref{PnmLxp1sa2}), (\ref{PnmLxp1la}),
(\ref{EstRnm}), (\ref{estsnmm0}) imply with $L=x+1$ and $\mu=1$ the estimate%
\begin{equation}
\left\vert \tilde{\omega}_{0}^{\left(  m\right)  }\left(  x\right)
\right\vert \leq\sqrt{3}\dfrac{\left(  m-1\right)  !}{x^{m}}M_{x,x+1}%
^{0,0,m-1}+\frac{M_{x,x+1}^{0,m,0}}{x^{m}}\log4+\frac{M_{x+1,\infty}^{0,m,0}%
}{x^{m}}+\frac{\left(  m-1\right)  !}{x^{m}}. \label{gnfin2}%
\end{equation}
The combination of (\ref{gnfin2}) with (\ref{estgntilde}) leads to%
\[
\left\vert \tilde{\omega}_{0}^{\left(  m\right)  }\left(  x\right)
\right\vert \leq m!\left(  \frac{\gamma}{x}\right)  ^{m}.
\]
This implies (\ref{estomegatildenm1}) and (\ref{estomegatildenm1la}) for
$m\geq1$ and $n=0$.

For $m\geq1$ and $n\geq1$, the choices $L=x+1$, and $\mu=1$ allow to estimate
the constant $K_{n,x+1}^{m,0}$ in (\ref{DefKL0}) by $K_{n,x+1}^{m,0}\leq
M_{x+1,x+1}^{1,m-1,n-1}$ and, in turn, we have%
\[
\left\vert S_{n,x+1,1}^{\left(  m\right)  }\left(  x\right)  \right\vert
\leq\frac{M_{x+1,x+1}^{1,m-1,n-1}}{nx^{m}}.
\]
Thus, the estimates (\ref{PnmLxp1la}), (\ref{PnmLxp1sa4}), (\ref{EstRnm})
imply
\begin{equation}
\left\vert \tilde{\omega}_{n}^{\left(  m\right)  }\left(  x\right)
\right\vert \leq\frac{M_{x,x+1}^{0,m,n}}{x^{m}}\log4+\sqrt{3}\frac
{M_{x,x+1}^{0,m-1,n-1}+M_{x,x+1}^{0,m-1,n}}{x^{m}}+\frac{M_{x+1,\infty
}^{0,m,n}}{x^{m}}+\frac{M_{x+1,x+1}^{1,m-1,n-1}}{nx^{m}}. \label{gnfin3}%
\end{equation}
Estimate (\ref{estgntilde}) allows to bound the terms $M_{a,b}^{\ell,m,n}$ in
(\ref{gnfin3}) which leads to (\ref{estomegatildenm1}) while the combination
with (\ref{estgntildela}) gives (\ref{estomegatildenm1la}).\medskip

\textbf{Case 3: Estimate for small argument}

For $m=0$, estimate (\ref{estomegatildenm1}) directly implies
(\ref{estomegatildenm1sa}). Note that the condition $0\leq m\leq\sqrt
{n}\left(  \frac{C-1}{C}\right)  $ implies $m=0$ for $n=0,1$. Hence, for the
following we assume $n\geq2$.

In the following, let $m\geq1$ and let $C$ and $\gamma$ be as explained in
statement (2a) of the theorem. We assume $1\leq m\leq\sqrt{n}\left(
\frac{C-1}{C}\right)  $ and restrict the range of $x$ to%
\begin{equation}
0<2x=:L\leq\frac{n}{C}\leq n-m\sqrt{n}. \label{recalll}%
\end{equation}
Then, (\ref{PnmL2x}) and (\ref{estgntildesa}) imply%
\begin{equation}
\left\vert P_{n,2x}^{\left(  m\right)  }\left(  x\right)  \right\vert
\leq\frac{M_{x,2x}^{0,m,n}}{x^{m}}\log4\leq\frac{m!}{\sqrt{n+1}}\left(
\frac{\gamma}{x}\right)  ^{m}, \label{finalsa1neu}%
\end{equation}
while (\ref{EstRnm}) and (\ref{estgntilde}) yield%
\[
\left\vert R_{n,2x,\mu}^{\left(  m\right)  }\left(  x\right)  \right\vert \leq
m!\left(  \frac{\gamma}{x}\right)  ^{m}\left(  n+1\right)  ^{\frac{m-1}{2}%
}\left(  \frac{1}{2}\right)  ^{2\mu}%
\]
By choosing $\mu:=\left\lfloor \frac{m\log\left(  n+1\right)  }{4\log
2}\right\rfloor $, we get%
\begin{equation}
\left\vert R_{n,2x,\mu}^{\left(  m\right)  }\left(  x\right)  \right\vert
\leq\frac{m!}{\sqrt{n+1}}\left(  \frac{\gamma}{x}\right)  ^{m}.
\label{finalsa2}%
\end{equation}
It remains to estimate the term $S_{n,L,\mu}^{\left(  m\right)  }\left(
x\right)  $. An estimate for the constant $K_{n,2x}^{m,r}$ for $0\leq r\leq
\mu-1$ follows from (\ref{DefKL0}), (\ref{estgntildesa}), and, by using
$2x=L\leq\frac{n}{C}$, via%
\begin{align}
K_{n,2x}^{m,r}  &  =L\sum_{\ell=0}^{\min\left(  m,2r\right)  }\tbinom{2r}%
{\ell}\frac{m!}{\left(  m-\ell\right)  !}\left\vert L^{m-\ell-1}%
g_{n-1-2r}^{\left(  m-1-\ell\right)  }\left(  L\right)  \right\vert
\label{Knest1}\\
&  \leq m!L\sum_{\ell=0}^{\min\left(  m,2r\right)  }\left(  m-\ell+1\right)
\tbinom{2r}{\ell}2^{m-\ell}\frac{1}{\sqrt{n-2r}}\left(  \frac{n-1-2r}%
{n-1-2r-L}\right)  ^{m-1-\ell}. \label{Knest1.5}%
\end{align}
Now, it holds
\[
\frac{n-1-2r}{2}\geq\frac{n+1-2\mu}{2}=\frac{1}{2}\left(  n-1-\left(
2\left\lfloor \frac{m\log\left(  n+1\right)  }{4\log2}\right\rfloor -2\right)
\right)  \overset{\text{(\ref{addcondm})}}{\geq}\frac{n-1}{4}.
\]
so that, by choosing $\gamma\geq6$ in (\ref{condxsmallargmts}), we obtain%
\begin{equation}
\frac{n-1-2r}{2}\geq\frac{n-1}{4}\overset{n\geq2}{\geq}\frac{n+1}{12}%
\geq2x=L\quad\text{so that (cf. (\ref{Knest1.5}))\quad}\frac{n-1-2r}%
{n-1-2r-L}\leq2. \label{Knest2}%
\end{equation}
The combination of (\ref{Knest1}) and (\ref{Knest2}) with $n-2r\geq n+2-2\mu$
and $L\leq n/C$ leads to%
\begin{align*}
K_{n,2x}^{m,r}  &  \leq\frac{\left(  m+1\right)  !n}{2C\sqrt{n+2-2\mu}}%
\sum_{\ell=0}^{\min\left(  m,2r\right)  }\tbinom{2r}{\ell}4^{m-\ell}\leq
\frac{\left(  m+1\right)  !n}{2C\sqrt{n+2-2\mu}}4^{m-2r}\sum_{\ell=0}%
^{2r}\tbinom{2r}{\ell}4^{2r-\ell}\\
&  =\frac{\left(  m+1\right)  !n}{2C\sqrt{n+2-2\mu}}4^{m-2r}5^{2r}\leq
\frac{\left(  n+1\right)  m!}{\sqrt{n+2-2\mu}}\gamma^{m+2r}%
\end{align*}
with a properly adjusted $\gamma$. The combination with (\ref{estsnm0}) and
$\left(  \frac{1}{2}\right)  ^{2r}\tbinom{2r}{r}\leq1$ leads to%
\[
\left\vert S_{n,L,\mu}^{\left(  m\right)  }\left(  x\right)  \right\vert \leq
m!\left(  \frac{\gamma}{x}\right)  ^{m}\frac{\left(  n+1\right)  }{\left(
n+2-2\mu\right)  ^{3/2}}\sum_{r=0}^{\mu-1}\left(  \frac{\gamma x}%
{n-2r}\right)  ^{2r}.
\]
Next, we use the additional condition (cf. (\ref{addcondm})) on $m$ such that%
\[
2r\leq2\mu-2=2\left\lfloor \frac{m\log\left(  n+1\right)  }{4\log
2}\right\rfloor -2\leq\frac{n-1}{2}%
\]
holds. This leads to
\[
\left\vert S_{n,L,\mu}^{\left(  m\right)  }\left(  x\right)  \right\vert \leq
m!\left(  \frac{\gamma}{x}\right)  ^{m}\frac{2^{3/2}}{\sqrt{n+1}}\sum
_{r=0}^{\mu-1}\left(  \frac{2\gamma x}{n+1}\right)  ^{2r}.
\]
By using the assumption $x\leq\frac{n+1}{4\gamma}$ as stated in (2a) of the
theorem we end up with the estimate%
\begin{equation}
\left\vert S_{n,L,\mu}^{\left(  m\right)  }\left(  x\right)  \right\vert \leq
m!\left(  \frac{\gamma}{x}\right)  ^{m}\frac{1}{\sqrt{n+1}}, \label{finalsa3}%
\end{equation}
again with a properly adjusted $\gamma$. The combination of (\ref{finalsa1neu}%
), (\ref{finalsa2}), and (\ref{finalsa3}) finally leads to the estimate
(\ref{estomegatildenm1sa}) for small argument.%
%TCIMACRO{\TeXButton{End Proof}{\endproof}}%
%BeginExpansion
\endproof
%EndExpansion

\section{Functional Estimates for Derivatives of $g_{n}\left(  t\right)
=\frac{t^{n}\operatorname*{e}^{-t}}{n!}$}

\begin{proposition}
\label{PropEstgnm}

\begin{enumerate}
\item General estimate. For $n\geq0$, $m\geq0$ and $\ell=0,1$, it holds for
all $t\geq0$%
\begin{equation}
\left\vert t^{m+\ell}g_{n}^{\left(  m\right)  }\left(  t\right)  \right\vert
\leq C_{m}\left(  n+1\right)  ^{\frac{m-1}{2}+\ell}\quad\text{with\quad}%
C_{m}:=\left(  4\operatorname*{e}\right)  ^{m+3}\left(  m+2\right)  !.
\label{estgntilde}%
\end{equation}

\item Refined estimates for small and large arguments.

\begin{enumerate}
\item Small argument. For $0\leq t\leq n-m\sqrt{n}$, the refined estimate
holds\footnote{For $n=0$, the condition $0\leq t\leq n-m\sqrt{n}$ implies
$t=m=0$ and the factor $\left(  \frac{n}{n-t}\right)  ^{m}$ is defined as $1$
for this case.}%
\begin{equation}
\left\vert t^{m}g_{n}^{\left(  m\right)  }\left(  t\right)  \right\vert
\leq2^{m+1}\frac{\left(  m+2\right)  !}{\sqrt{n+1}}\left(  \frac{n}%
{n-t}\right)  ^{m}. \label{estgntildesa}%
\end{equation}

\item Large argument. For $n\geq0$ and $m\geq2\log\left(  n+1\right)  $, the
refined estimate%
\begin{equation}
\left\vert t^{m+\ell}g_{n}^{\left(  m\right)  }\left(  t\right)  \right\vert
\leq4\left(  \frac{3}{\ln\frac{1}{c}}\right)  ^{m}\left(  m+2\right)  !\left(
n+1\right)  ^{\ell}\qquad\forall t\geq\left\{
\begin{array}
[c]{ll}%
0 & n=0,1\\
n+m\left(  \sqrt{n}+2\right)  & n\geq2
\end{array}
\right.  \label{estgntildela}%
\end{equation}
holds with $c:=9/10.$
\end{enumerate}

\item Exponential decay. For $m=0$ and $t\geq n+\sqrt{n}$, we obtain%
\begin{equation}
g_{n}\left(  t\right)  \leq\frac{1}{\sqrt{n+1}}\exp\left(  \sqrt{n}-\frac
{t}{1+\sqrt{n}}\right)  . \label{estgnasymptotic}%
\end{equation}

\end{enumerate}
\end{proposition}%

%TCIMACRO{\TeXButton{Proof}{\proof}}%
%BeginExpansion
\proof
%EndExpansion
We start to prove some special cases.

\textbf{Case }$n=0,1$.

For $n=0$, and $\ell=0,1$, it holds%
\begin{equation}
\left\vert t^{m+\ell}g_{0}^{\left(  m\right)  }\left(  t\right)  \right\vert
=t^{m+\ell}\operatorname*{e}\nolimits^{-t}\leq\frac{\left(  m+\ell\right)
!}{\sqrt{m+\ell+1}}\qquad\forall t\geq0, \label{estg0m}%
\end{equation}
so that (\ref{estgntilde}) holds for all $C_{m}\geq\frac{\left(
m+\ell\right)  !}{\sqrt{m+\ell+1}}$. This also implies (\ref{estgntildesa})
and (\ref{estgnasymptotic}) for $n=0$.

For $n=1$, we get%
\[
\left(  -t\right)  ^{m+\ell}g_{1}^{\left(  m\right)  }\left(  t\right)
=\left(  -t\right)  ^{m+\ell}\left(  t\operatorname*{e}\nolimits^{-t}\right)
^{\left(  m\right)  }=\left(  -1\right)  ^{\ell}\left(  t-m\right)  t^{m+\ell
}\operatorname*{e}\nolimits^{-t}.
\]
By estimating $t^{m+\ell+1}\operatorname*{e}^{-t}$ and $t^{m+\ell
}\operatorname*{e}^{-t}$ as in (\ref{estg0m}) we get%
\begin{equation}
\left\vert t^{m+\ell}g_{1}^{\left(  m\right)  }\left(  t\right)  \right\vert
\leq2\frac{\left(  m+\ell+1\right)  !}{\sqrt{m+\ell+1}} \label{estg0m1}%
\end{equation}
so that (\ref{estgntilde}) holds for $n=1$ if $C_{m}\geq2^{\frac{3-m}{2}-\ell
}\frac{\left(  m+\ell+1\right)  !}{\sqrt{m+\ell+1}}$. Note that this also
implies (\ref{estgntildesa}) and (\ref{estgnasymptotic}) for $n=1$.

For the rest of the proof we assume $n\geq2$. Note that (\ref{defDmsa})
directly implies (\ref{estgntildesa}) for $n\geq2$ and it remains to prove the
remaining inequalities.\medskip

\textbf{Case }$m=0,1$.

For $m=0$ and $\ell=0,1$ the function $t^{\ell}g_{n}\left(  t\right)  $ has
its extremum at $t=n+\ell$, i.e.,%
\[
\left\vert t^{\ell}g_{n}\left(  t\right)  \right\vert \leq\frac{\left(
n+\ell\right)  ^{n+\ell}\operatorname*{e}^{-n-\ell}}{n!}\qquad\forall t\geq0.
\]
Since $n\geq1$, Stirling's formula gives us%
\[
\left\vert t^{\ell}g_{n}\left(  t\right)  \right\vert \leq\frac{\left(
n+1\right)  ^{\ell}}{\sqrt{n+\ell+1}}\leq\left(  n+1\right)  ^{\ell-1/2}%
\]
so that (\ref{estgntilde}) holds for this case if $C_{0}\geq1$.

For $m=1$ and $\ell=0,1$, the function $t^{1+\ell}g_{n}^{\left(  1\right)
}\left(  t\right)  =\left(  n-t\right)  \frac{t^{n+\ell}\operatorname*{e}%
^{-t}}{n!}$ has its extrema at $t_{\pm}=\left(  n+\frac{\ell+1}{2}\right)
\left(  1\pm\delta_{n,\ell}\right)  $ with $\delta_{n,\ell}=\frac
{\sqrt{n+\frac{\left(  \ell+1\right)  ^{2}}{4}}}{n+\frac{\ell+1}{2}}$. Hence,
with Stirling's formula we obtain%
\begin{align}
\left\vert t_{\pm}^{1+\ell}g_{n}^{\prime}\left(  t_{\pm}\right)  \right\vert
&  =\left(  \sqrt{n+\frac{\left(  \ell+1\right)  ^{2}}{4}}\pm\frac{\ell+1}%
{2}\right)  \frac{\left(  \left(  n+\frac{\ell+1}{2}\right)  \left(
1\pm\delta_{n,\ell}\right)  \right)  ^{n+\ell}\operatorname*{e}^{-\left(
n+\frac{\ell+1}{2}\right)  \left(  1\pm\delta_{n,\ell}\right)  }}%
{n!}\nonumber\\
&  \leq\frac{1}{\sqrt{2\pi}}\left(  \sqrt{1+\frac{\left(  \ell+1\right)  ^{2}%
}{4n}}\pm\frac{\ell+1}{2\sqrt{n}}\right)  \left(  n+\frac{\ell+1}{2}\right)
^{\ell}\left(  1+\frac{\ell+1}{2n}\right)  ^{n}\times\label{estg1t}\\
&  \times\left(  \left(  \left(  1\pm\delta_{n,\ell}\right)  \right)
^{n+\ell}\operatorname*{e}\nolimits^{-\left(  \frac{\ell+1}{2}\pm\left(
n+\frac{\ell+1}{2}\right)  \delta_{n,\ell}\right)  }\right)  . \label{35}%
\end{align}
Since $n\geq1$ and $\ell=0,1$, we get%
\[
\sqrt{1+\frac{\left(  \ell+1\right)  ^{2}}{4n}}\pm\frac{\ell+1}{2\sqrt{n}}%
\leq\sqrt{2}+1,\text{\quad}\left(  1+\frac{\ell+1}{2n}\right)  ^{n}%
\leq\operatorname*{e}\text{ and\ }\left(  n+\frac{\ell+1}{2}\right)  ^{\ell
}=\left(  n+1\right)  ^{\ell}.
\]
The last factor in (\ref{estg1t}) is considered first with \textquotedblleft%
+\textquotedblright\ signs and can be estimated by%
\begin{align*}
\left(  \left(  1+\delta_{n,\ell}\right)  \right)  ^{n+\ell}\operatorname*{e}%
\nolimits^{-\left(  \frac{\ell+1}{2}+\left(  n+\frac{\ell+1}{2}\right)
\delta_{n,\ell}\right)  }  &  =\operatorname*{e}\nolimits^{\left(
n+\ell\right)  \log\left(  1+\delta_{n}\right)  }\operatorname*{e}%
\nolimits^{-\left(  \frac{\ell+1}{2}+\left(  n+\frac{\ell+1}{2}\right)
\delta_{n,\ell}\right)  }\\
&  \leq\operatorname*{e}\nolimits^{\left(  n+\ell\right)  \delta_{n,\ell
}-\left(  \frac{\ell+1}{2}+\left(  n+\frac{\ell+1}{2}\right)  \delta_{n,\ell
}\right)  }=\operatorname*{e}\nolimits^{-c_{n,\ell}}%
\end{align*}
with%
\[
c_{n,\ell}=\frac{1-\ell}{2}\delta_{n,\ell}+\frac{\left(  \ell+1\right)  }%
{2}\geq0
\]
so that, in this case, the last factor in (\ref{estg1t}), (\ref{35}) is
bounded by $1$. For the \textquotedblleft$-$\textquotedblright\ signs, we get%
\[
\left(  n+\ell\right)  \log\left(  1-\delta_{n,\ell}\right)  =-\left(
n+\ell\right)  \log\left(  1+\frac{\delta_{n,\ell}}{1-\delta_{n,\ell}}\right)
\leq-\left(  n+\ell\right)  \delta_{n,\ell}%
\]
so that%
\begin{align*}
\left(  \left(  1-\delta_{n,\ell}\right)  \right)  ^{n+\ell}\operatorname*{e}%
\nolimits^{-\left(  \frac{\ell+1}{2}-\left(  n+\frac{\ell+1}{2}\right)
\delta_{n,\ell}\right)  }  &  \leq\exp\left(  -\left(  n+\ell\right)
\delta_{n,\ell}-\left(  \frac{\ell+1}{2}-\left(  n+\frac{\ell+1}{2}\right)
\delta_{n,\ell}\right)  \right) \\
&  =\exp\left(  \frac{1-\ell}{2}\delta_{n,\ell}-\frac{\ell+1}{2}\right)
\end{align*}
Since $\delta_{n,\ell}\leq\frac{1}{\sqrt{2}}$ we arrive at the estimate
\[
\left\vert t^{1+\ell}g_{n}^{\prime}\left(  t\right)  \right\vert \leq
\frac{3\operatorname*{e}}{\sqrt{2\pi}}\left(  n+1\right)  ^{\ell}%
\]
and this proves (\ref{estgntilde}) for $m=1$.\medskip

\textbf{Case }$0\leq t\leq n-\sqrt{n}$.

Proposition \ref{Propgm2} implies for $\ell=0,1$ and $t\leq n$
\begin{equation}
\left\vert t^{m+\ell}g_{n}^{\left(  m\right)  }\left(  t\right)  \right\vert
\leq\left(  m+2\right)  !\operatorname*{e}\nolimits^{m}n^{\frac{m-1}{2}+\ell}.
\label{case0tsqrtn}%
\end{equation}
\medskip

\textbf{Case }$n-\sqrt{n}\leq t\leq n+\sqrt{n}$.

We start with the simple recursions%

\[
g_{n}^{\prime}=g_{n-1}-g_{n}\quad\text{and\quad}g_{n}=\frac{t}{n}g_{n-1}%
\]
from which we conclude $g_{n}^{\prime}=\left(  1-\frac{t}{n}\right)  g_{n-1}$.
By differentiating this relation $m$ times we get via Leibniz' rule%
\[
g_{n}^{\left(  m+1\right)  }=\left(  1-\frac{t}{n}\right)  g_{n-1}^{\left(
m\right)  }-\frac{m}{n}g_{n-1}^{\left(  m-1\right)  }\quad\text{with (cf.
(\ref{estg0m})) \quad}\left\vert g_{0}^{\left(  m\right)  }\right\vert
\leq1,\quad\left\vert g_{1}^{\left(  m\right)  }\right\vert \leq m.
\]
For $n=0,1$ we get%
\[
\left\vert g_{0}^{\left(  m\right)  }\left(  t\right)  \right\vert
=\operatorname*{e}\nolimits^{-t}\leq1\quad\text{and\quad}g_{1}^{\left(
m\right)  }\left(  t\right)  =\left\vert \left(  t-m\right)  \exp\left(
-t\right)  \right\vert \leq m.
\]
It is easy to verify that the coefficients $A_{n}^{\left(  m\right)  }$ in the
recursion%
\begin{equation}
A_{n}^{\left(  m\right)  }:=\left\{
\begin{array}
[c]{lc}%
1 & n=0,m\in\mathbb{N}_{0},\\
m & n=1,m\in\mathbb{N}_{0},\\
\dfrac{1}{\sqrt{n+1}} & m=0,n\in\mathbb{N}_{0},\\
\frac{\sqrt{2}}{n+1} & m=1,n\in\mathbb{N}_{\geq2},\\
A_{n}^{\left(  m\right)  }=\frac{1}{\sqrt{n}}A_{n-1}^{\left(  m-1\right)
}+\frac{m-1}{n}A_{n-1}^{\left(  m-2\right)  } & n\geq2,m\geq2
\end{array}
\right.  \label{defanmrec}%
\end{equation}
majorate $\left\vert g_{n}^{\left(  m\right)  }\right\vert $. For the estimate
of $A_{n}^{\left(  m\right)  }$ we distinguish between two cases. Recall that
we may restrict to the cases $m\geq2$ and $n\geq2$.

\textbf{a) }Let $n\geq m/2$. Then it holds%
\begin{equation}
A_{n}^{\left(  m\right)  }\leq\frac{m!a^{m}}{n^{\left(  m+1\right)  /2}}%
\quad\text{for any\quad}a\geq1+\sqrt{3}\text{.} \label{estimateAnm}%
\end{equation}
This is proved by induction: It is easy to see that the right-hand side in
(\ref{estimateAnm}) majorates $A_{n}^{\left(  m\right)  }$ for the first four
cases in (\ref{defanmrec}). Then, by induction we get%
\begin{align}
\frac{1}{\sqrt{n}}A_{n-1}^{\left(  m-1\right)  }+\frac{m-1}{n}A_{n-1}^{\left(
m-2\right)  }  &  \leq\frac{1}{\sqrt{n}}\frac{\left(  m-1\right)  !a^{m-1}%
}{\left(  n-1\right)  ^{m/2}}+\frac{m-1}{n}\frac{\left(  m-2\right)  !a^{m-2}%
}{\left(  n-1\right)  ^{\left(  m-1\right)  /2}}\nonumber\\
&  \leq\frac{m!a^{m}}{n^{\left(  m+1\right)  /2}}\left(  \frac{1}{ma^{2}%
}\left(  a\left(  \frac{n}{n-1}\right)  ^{m/2}+\left(  \frac{n}{n-1}\right)
^{\left(  m-1\right)  /2}\right)  \right)  . \label{estAnm2}%
\end{align}
For $n\geq m/2$, we have%
\[
\left(  \frac{n}{n-1}\right)  ^{\left(  m-1\right)  /2}\leq\left(  \frac
{n}{n-1}\right)  ^{m/2}\leq\left(  \frac{n}{n-1}\right)  ^{n}\leq4
\]
so that the factor $\left(  \ldots\right)  $ in the right-hand side of
(\ref{estAnm2}) can be estimated from above by $2\left(  a+1\right)  /a^{2}$,
which is $\leq1$ for $a\geq1+\sqrt{3}$. Thus, estimate (\ref{estimateAnm}) is
proved .

\textbf{b) }For $n<m/2$, it holds%
\begin{equation}
A_{n}^{\left(  m\right)  }\leq G_{n}^{\left(  m\right)  }:=\frac{2^{n}\left(
m-1\right)  !!}{n!\left(  m-1-2n\right)  !!}. \label{defgnms}%
\end{equation}
This can be seen by first observing that this holds for $n=0$. Next we prove
the auxiliary statement\textbf{: }For $1\leq n<m/2$, it holds%
\begin{equation}
\frac{1}{\sqrt{n}}G_{n-1}^{\left(  m\right)  }\leq\frac{m}{n}G_{n-1}^{\left(
m-1\right)  }. \label{auxstatpos}%
\end{equation}
This is equivalent to%
\[
Q_{n}^{\left(  m\right)  }:=\frac{\sqrt{n}}{m}\frac{\left(  m-1\right)
!!}{\left(  m-2\right)  !!}\frac{\left(  m-2n\right)  !!}{\left(
m+1-2n\right)  !!}\leq1.
\]
It is easy to see that the quotient $\frac{\left(  m-2n\right)  !!}{\left(
m+1-2n\right)  !!}$ increases with increasing $1\leq n<m/2$ so that
$Q_{n}^{\left(  m\right)  }$ can be bounded from above by setting
$n=\frac{m-1}{2}$:%
\[
Q_{n}^{\left(  m\right)  }\leq\frac{\sqrt{m-1}}{2\sqrt{2}m}\frac{\left(
m-1\right)  !!}{\left(  m-2\right)  !!}\overset{\text{Lem.
\ref{LemDoubleFactorial}}}{\leq}\frac{m-1}{\sqrt{2}m}\leq1
\]
and the auxiliary statement is proved.\medskip

Hence, by induction it holds%
\begin{align*}
A_{n}^{\left(  m+1\right)  }  &  =\frac{1}{\sqrt{n}}A_{n-1}^{\left(  m\right)
}+\frac{m}{n}A_{n-1}^{\left(  m-1\right)  }\leq\frac{1}{\sqrt{n}}%
G_{n-1}^{\left(  m\right)  }+\frac{m}{n}G_{n-1}^{\left(  m-1\right)  }%
=2\frac{m}{n}G_{n-1}^{\left(  m-1\right)  }\\
&  =\frac{2^{n}m!!}{n!\left(  m-2n\right)  !!}=G_{n}^{\left(  m+1\right)  }%
\end{align*}
and (\ref{defgnms}) is proved. \medskip

\textbf{c) }For $1\leq n<m/2$, we will show that%
\begin{equation}
G_{n}^{\left(  m\right)  }\leq\frac{m!a^{m}}{n^{\left(  m+1\right)  /2}}
\label{estgnmfinal}%
\end{equation}
holds. For the smallest values of $m$, i.e., $m=2n+1$, this follows from
Stirling's formula with $a\geq1+\sqrt{3}$%
\begin{align*}
G_{n}^{\left(  2n+1\right)  }  &  =4^{n}\leq\left(  \frac{2a}%
{\operatorname*{e}}\right)  ^{2n}\leq2\sqrt{\pi}\left(  \frac{2a}%
{\operatorname*{e}}\right)  ^{2n+1}\left(  n+\frac{1}{2}\right)
^{n+1/2}\left(  1+\frac{1}{2n}\right)  ^{n+1}\\
&  =a^{2n+1}\sqrt{2\pi}\frac{\left(  2n+1\right)  ^{2n+3/2}\operatorname*{e}%
^{-\left(  2n+1\right)  }}{n^{n+1}}\leq\frac{a^{m}m!}{n^{\left(  m+1\right)
/2}}.
\end{align*}
Hence, for $m>2n+1$ we get by induction%
\[
G_{n}^{\left(  m\right)  }\overset{\text{(\ref{auxstatpos})}}{\leq}\frac
{m}{\sqrt{n+1}}G_{n}^{\left(  m-1\right)  }\overset{\text{induction}}{\leq
}\frac{m}{\sqrt{n+1}}\frac{\left(  m-1\right)  !a^{m-1}}{n^{m/2}}\leq
\frac{m!a^{m}}{n^{\left(  m+1\right)  /2}}%
\]
and (\ref{estgnmfinal}) is proved.

Since $n\geq$ $2$, it follows%
\[
\left\vert g_{n}^{\left(  m\right)  }\left(  t\right)  \right\vert \leq
\frac{\delta_{m}}{\left(  n+1\right)  ^{(m+1)/2}}\quad\text{with\quad}%
\delta_{m}=m!\left(  1+\sqrt{3}\right)  ^{m}\left(  \frac{3}{2}\right)
^{\left(  m+1\right)  /2}.
\]
From $t^{m+\ell}\leq\left(  n+\sqrt{n}\right)  ^{m+\ell}\leq2^{m+\ell}\left(
n+1\right)  ^{m+\ell}$, the assertion follows:%
\[
\left\vert t^{m+\ell}g_{n}^{\left(  m\right)  }\left(  t\right)  \right\vert
\leq2^{m+1}\delta_{m}\left(  n+1\right)  ^{\frac{m-1}{2}+\ell}.
\]

\textbf{Case }$t\geq n+\sqrt{n}$.

Estimate (\ref{estgntilde}) in this case follows from Proposition
\ref{Propgnlargt}.\medskip

The estimate (\ref{estgnasymptotic}) follows trivially, for $n=0$, from the
definition of $g_{n}$ and, for $n\geq2$, directly from Proposition
\ref{Propgnlargt}. For $n=1$, the estimate follows by observing that
(\ref{gnt}) also holds for $n=1$. The refined estimate (\ref{estgntildela})
follows for $n\geq2$ from (\ref{tmgplusla}) and the case $n=0,1$ have been
treated already at the beginning of the proof.%
%TCIMACRO{\TeXButton{End Proof}{\endproof}}%
%BeginExpansion
\endproof
%EndExpansion

\begin{proposition}
\label{Propgm2}Let $n\geq2$ and $m\geq0$. For $0\leq t\leq n-\sqrt{n}$, it
holds%
\begin{equation}
\left\vert t^{m}g_{n}^{\left(  m\right)  }\left(  t\right)  \right\vert
\leq\operatorname*{e}\nolimits^{m}\left(  m+2\right)  !n^{\frac{m-1}{2}}.
\label{defDm}%
\end{equation}
For $0\leq t\leq n-m\sqrt{n}$, the refined estimate holds%
\begin{equation}
\left\vert t^{m}g_{n}^{\left(  m\right)  }\left(  t\right)  \right\vert
\leq2^{m+1}\frac{\left(  m+2\right)  !}{\sqrt{n+1}}\left(  \frac{n}%
{n-t}\right)  ^{m}. \label{defDmsa}%
\end{equation}

\end{proposition}%

%TCIMACRO{\TeXButton{Proof}{\proof}}%
%BeginExpansion
\proof
%EndExpansion
Note that%
\begin{equation}
t^{m}g_{n}^{\left(  m\right)  }\left(  t\right)  =g_{n}\left(  t\right)
s_{n,m}\left(  t\right)  \label{repgknm1}%
\end{equation}
with%
\begin{equation}
s_{n,m}\left(  t\right)  =\sum_{\ell=0}^{\min\left(  m,n\right)  }\tbinom
{m}{\ell}\left(  -1\right)  ^{m-\ell}\frac{n!t^{m-\ell}}{\left(
n-\ell\right)  !}. \label{defsnm}%
\end{equation}

\textbf{Estimate of }$g_{n}$.

We set $t=n/c$ for some $c=1+\frac{\delta}{n-\delta}$ and $0\leq\delta<n.$
Stirling's formula gives us%
\[
g_{n}\left(  \frac{n}{c}\right)  \leq\frac{w^{n}\left(  c\right)  }{\sqrt
{n+1}},
\]
where%
\begin{align*}
w^{n}\left(  c\right)   &  :=\left(  \frac{\exp\left(  1-\frac{1}{c}\right)
}{c}\right)  ^{n}=\left(  1-\frac{\delta}{n}\right)  ^{n}\exp\left(
\delta\right) \\
&  =\exp\left(  n\log\left(  1-\frac{\delta}{n}\right)  +\delta\right)
=\exp\left(  -n\sum_{k=2}^{\infty}\left(  \frac{\delta}{n}\right)
^{k}/k\right)  \leq\exp\left(  -\frac{\delta^{2}}{2n}\right)  .
\end{align*}
Note that the range $t\in\left[  0,n-\sqrt{n}\right]  $ corresponds to the
range of $\delta\in\left[  \sqrt{n},n\right]  $. Thus%
\begin{equation}
g_{n}\left(  n-\delta\right)  \leq\frac{\exp\left(  -\delta^{2}/\left(
2n\right)  \right)  }{\sqrt{n+1}}\quad\forall\delta\in\left[  \sqrt
{n},n\right]  . \label{repgknm2}%
\end{equation}
This proves the assertion (\ref{defDm}) for $m=0$ so that, for the following,
we assume $m\geq1$.

\textbf{Estimate of }$s_{n,m}$ \textbf{for }$m\geq1$.

The definition of $s_{n,m}$ directly leads to the estimate%
\begin{equation}
\left\vert s_{n,m}\left(  t\right)  \right\vert \leq\sum_{\ell=0}^{m}%
\tbinom{m}{\ell}n^{\ell}t^{m-\ell}\leq\left(  t+n\right)  ^{m}\text{.}
\label{genestsnm}%
\end{equation}

Next this estimate will be refined for $0\leq t\leq n-\sqrt{n}$. We set
$\delta:=n-t$ and introduce the function $\tilde{s}_{n,m}\left(
\delta\right)  =\left(  -1\right)  ^{m}s_{n,m}\left(  n-\delta\right)  $ so
that an estimate of $\left\vert \tilde{s}_{n,m}\right\vert $ at $\delta$
implies the same estimate of $\left\vert s_{n,m}\right\vert $ at $n-\delta$.
For later use, we will estimate $\tilde{s}_{n,m}\left(  \delta\right)  $ not
only for $\delta\in\left[  \sqrt{n},n\right]  $ but for all $\delta
\in\mathbb{R}$ with $\left\vert \delta\right\vert \geq\sqrt{n}$.

From (\ref{defsnm}) one concludes that $\tilde{s}_{n,m}$ satisfies the
recursion%
\begin{equation}
\tilde{s}_{n,m+1}\left(  \delta\right)  =\left(  m-\delta\right)  \tilde
{s}_{n,m}\left(  t\right)  +\left(  n-\delta\right)  \tilde{s}_{n,m}^{\prime
}\left(  \delta\right)  \quad\text{with\quad}\tilde{s}_{0,0}:=1. \label{reks}%
\end{equation}
By inspection of (\ref{defsnm}) we conclude that%
\begin{equation}
\tilde{s}_{n,m}\left(  \delta\right)  =\sum_{\ell=0}^{m}n^{\ell}p_{\ell
,m}\left(  \delta\right)  \mathbf{,} \label{defstilde}%
\end{equation}
where $p_{\ell,m}\in\mathbb{P}_{m-\ell}$. From (\ref{reks}) we obtain the
recursion%
\begin{equation}
p_{\ell,m+1}\left(  \delta\right)  =\left(  m-\delta\right)  p_{\ell,m}\left(
\delta\right)  -\delta p_{\ell,m}^{\prime}\left(  \delta\right)  +p_{\ell
-1,m}^{\prime}\left(  \delta\right)  \quad\text{with\quad}p_{0,0}:=1,
\label{rek3}%
\end{equation}
where we formally set $p_{-1,m}=0$ and $p_{\ell,m}=0$ for $\ell>m$. It is easy
to prove by induction that%
\[
p_{0,m}\left(  \delta\right)  =\left(  -1\right)  ^{m}\delta^{m}%
\]
and $p_{\ell,m}\in\mathbb{P}_{m-2\ell}$, where $\mathbb{P}_{\ell}:=\left\{
0\right\}  $ for $\ell<0$. Hence,%
\begin{equation}
\tilde{s}_{n,m}\left(  \delta\right)  =\sum_{\ell=0}^{\left\lfloor
m/2\right\rfloor }n^{\ell}p_{\ell,m}\left(  \delta\right)  . \label{defplm}%
\end{equation}
Next, we will estimate the polynomial $p_{\ell,m}$. We write%
\[
p_{\ell,m}\left(  \delta\right)  =\sum_{k=0}^{m-2\ell}a_{\ell,m,k}\delta
^{k}\quad\text{with\quad}a_{\ell,m,k}=c_{\ell,m-k}(-1)^{k+\ell}\frac{m!}{k!}.
\]
Plugging this ansatz into (\ref{rek3}) gives%
\[
a_{\ell,m+1,k}+a_{\ell,m,k-1}=\left(  m-k\right)  a_{\ell,m,k}+\left(
k+1\right)  a_{\ell-1,m,k+1}%
\]
and, in turn, the recursion\footnote{Particular results are%
\begin{align*}
c_{0,k}  &  =\delta_{0,k}\\
c_{1,k}  &  =\frac{1}{k}\\
c_{2,k}  &  =\frac{H_{k-2}-1}{k}\\
c_{3,k}  &  =-\frac{-1+2H_{k-3}-H_{k-2}^{2}+H_{k-2,2}}{2k}%
\end{align*}
with the harmonic numbers $H_{n,r}=\sum_{\ell=1}^{n}1/\ell^{r}$.}%
\[
c_{\ell,k+1}=\frac{k}{k+1}c_{\ell,k}+\frac{c_{\ell-1,k-1}}{k+1}\quad
\text{with\quad}c_{\ell,2\ell}=\frac{1}{2^{\ell}\ell!}.
\]
By induction it is easy to prove that%
\[
c_{\ell,k}\leq1\quad\forall k\geq2\ell
\]
so that%
\[
\left\vert p_{\ell,m}\left(  \delta\right)  \right\vert \leq m!\sum
_{k=0}^{m-2\ell}\frac{\left\vert \delta\right\vert ^{k}}{k!}.
\]
Note that $t^{k-1}/\left(  k-1\right)  !\leq t^{r}/r!$ for all $t\geq r\geq
k\geq1$ so that%
\begin{equation}
\left\vert p_{\ell,m}\left(  \delta\right)  \right\vert \leq\left\{
\begin{array}
[c]{ll}%
\dfrac{\left(  m+1\right)  !}{\left(  m-2\ell\right)  !}\left\vert
\delta\right\vert ^{m-2\ell} & \left\vert \delta\right\vert \geq m-2\ell,\\
m!\operatorname*{e}^{\left\vert \delta\right\vert } & \delta\in\mathbb{R}.
\end{array}
\right.  \label{estplmabs}%
\end{equation}
Hence, it holds for $\left\vert \delta\right\vert \geq\sqrt{n}$%
\begin{align}
\left\vert s_{n,m}\left(  n-\delta\right)  \right\vert =\left\vert \tilde
{s}_{n,m}\left(  \delta\right)  \right\vert  &  \leq\left\{
\begin{array}
[c]{lc}%
n^{m/2}\left(  m+2\right)  !\max\limits_{0\leq\ell\leq m/2}\dfrac{\left(
\left\vert \delta\right\vert /\sqrt{n}\right)  ^{m-2\ell}}{\left(
m-2\ell\right)  !} & \left\vert \delta\right\vert \geq m,\\
n^{m/2}\left(  m+1\right)  !\operatorname*{e}^{\left\vert \delta\right\vert }
& \delta\in\mathbb{R},
\end{array}
\right. \nonumber\\
&  \leq\left\{
\begin{array}
[c]{ll}%
\left(  m+2\right)  ^{2}\delta^{m} & m\sqrt{n}\leq\delta,\\
n^{m/2}\left(  m+2\right)  !\operatorname*{e}^{\left\vert \delta\right\vert
/\sqrt{n}} & m\leq\delta\leq m\sqrt{n},\\
n^{m/2}\left(  m+1\right)  !\operatorname*{e}^{\left\vert \delta\right\vert }
& \delta\in\mathbb{R}.
\end{array}
\right.  \label{repgknm3}%
\end{align}

\textbf{Estimate of }$t^{m}g_{n}^{\left(  m\right)  }\left(  t\right)  .$

Let $t=n-\delta$ for some $\delta\in\left[  \sqrt{n},n\right]  $.

\textbf{Case }$\sqrt{n}\leq\delta\leq m$. The combination of the third case in
(\ref{repgknm3}) with (\ref{repgknm2}) yields%
\begin{equation}
\left\vert t^{m}g_{n}^{\left(  m\right)  }\left(  t\right)  \right\vert \leq
n^{m/2}\left(  m+1\right)  !\frac{\exp\left(  \delta-\frac{\delta^{2}}%
{2n}\right)  }{\sqrt{n+1}}. \label{esttmgnm}%
\end{equation}
From $\delta\leq m$ we conclude that%
\[
\left\vert t^{m}g_{n}^{\left(  m\right)  }\left(  t\right)  \right\vert \leq
n^{\frac{m-1}{2}}\left(  m+1\right)  !\operatorname*{e}\nolimits^{m}.
\]

\textbf{Case }$\max\left(  m,\sqrt{n}\right)  \leq\delta\leq\sqrt{n}%
\min\left(  m,\sqrt{n}\right)  $. Here we get from the second case in
(\ref{repgknm3}) and (\ref{repgknm2}) the estimate%
\[
\left\vert t^{m}g_{n}^{\left(  m\right)  }\left(  t\right)  \right\vert \leq
n^{m/2}\left(  m+2\right)  !\frac{\exp\left(  \frac{\delta}{\sqrt{n}}%
-\frac{\delta^{2}}{2n}\right)  }{\sqrt{n+1}}.
\]
The exponent is monotonously decreasing for $\delta\geq\sqrt{n}$ so that%
\[
\left\vert t^{m}g_{n}^{\left(  m\right)  }\left(  t\right)  \right\vert
\leq\sqrt{\operatorname*{e}}n^{\frac{m-1}{2}}\left(  m+2\right)  !.
\]

\textbf{Case }$m\sqrt{n}\leq\delta\leq n$. In this case, we obtain, by using
$\operatorname*{e}^{t}\geq t^{k}/k!$ for $t\geq0$, from the first case in
(\ref{repgknm3}) and (\ref{repgknm2})%
\[
\left\vert t^{m}g_{n}^{\left(  m\right)  }\left(  t\right)  \right\vert
\leq\frac{\left(  m+2\right)  ^{2}}{\sqrt{n+1}}\frac{\delta^{m}}{\exp\left(
\delta^{2}/\left(  2n\right)  \right)  }\leq2^{m}\frac{\left(  m+2\right)
^{2}}{\sqrt{n+1}}m!\left(  \frac{n}{\delta}\right)  ^{m}.
\]%
%TCIMACRO{\TeXButton{End Proof}{\endproof}}%
%BeginExpansion
\endproof
%EndExpansion

\begin{proposition}
\label{Propgnlargt}Let $n\geq2$.

\begin{enumerate}
\item For $t\geq n+\sqrt{n}$, it holds%
\begin{equation}
g_{n}\left(  t\right)  \leq\frac{1}{\sqrt{n+1}}\exp\left(  \sqrt{n}\left(
1-\frac{t}{n+\sqrt{n}}\right)  \right)  . \label{expdecay}%
\end{equation}

\item For $m\geq0$ and $\ell=0,1$, we get the estimate%
\begin{equation}
\left\vert t^{m+\ell}g_{n}^{\left(  m\right)  }\left(  t\right)  \right\vert
\leq\left(  4\operatorname*{e}\right)  ^{m+3}\left(  m+2\right)
!n^{\frac{m-1}{2}+\ell}\quad\forall t\geq n+\sqrt{n}. \label{tmgplus}%
\end{equation}

\item For $t\geq n+m\left(  \sqrt{n}+2\right)  $ and $m\geq2\log n$, the
refined estimate%
\begin{equation}
\left\vert t^{m+\ell}g_{n}^{\left(  m\right)  }\left(  t\right)  \right\vert
\leq\left(  \frac{3}{\ln\frac{1}{c}}\right)  ^{m}\left(  m+2\right)  !\left(
4n\right)  ^{\ell}. \label{tmgplusla}%
\end{equation}
holds with $c:=9/10.$
\end{enumerate}
\end{proposition}

%

%TCIMACRO{\TeXButton{Proof}{\proof}}%
%BeginExpansion
\proof
%EndExpansion
For any $0<c\leq1$, we write%
\[
g_{n}\left(  t\right)  =\frac{t^{n}\exp\left(  -ct\right)  }{n!}\exp\left(
-\left(  1-c\right)  t\right)  .
\]
The first fraction has its maximum at $t=n/c$ so that Stirling's formula leads
to%
\begin{equation}
g_{n}\left(  t\right)  \leq\frac{1}{\sqrt{n+1}}\left(  \frac{1}{c}\right)
^{n}\exp\left(  -\left(  1-c\right)  t\right)  . \label{gnt0}%
\end{equation}
We choose $c=\frac{n}{n+\sqrt{n}}$ and obtain%
\begin{equation}
g_{n}\left(  t\right)  \leq\frac{1}{\sqrt{n+1}}\left(  1+\frac{1}{\sqrt{n}%
}\right)  ^{n}\exp\left(  -\frac{\sqrt{n}}{n+\sqrt{n}}t\right)  \leq\frac
{1}{\sqrt{n+1}}\exp\left(  \sqrt{n}\left(  1-\frac{t}{n+\sqrt{n}}\right)
\right)  \label{gnt}%
\end{equation}
which shows the exponential decay for $t\geq n+\sqrt{n}$ (cf. (\ref{expdecay}%
)). This also implies (\ref{tmgplus}) for $m=\ell=0$. For $m=0$ and $\ell=1$,
we obtain%
\[
\left\vert tg_{n}\left(  t\right)  \right\vert \leq\frac{t^{n+1}%
\operatorname*{e}\nolimits^{-t}}{n!}\leq\frac{\left(  n+1\right)
^{n+1}\operatorname*{e}^{-n-1}}{n!}\overset{\text{Stirling}}{\leq}\sqrt{n+1}%
\]
so that (\ref{tmgplus}) is satisfied for $m=0$, $\ell=0,1$.\medskip

For the rest of the proof, we assume $m\geq1$ and employ the representation
(\ref{repgknm1}).

\textbf{Estimate of }$t^{m+\ell}g_{n}^{\left(  m\right)  }\left(  t\right)
$\textbf{.}

The combination of (\ref{repgknm1}), (\ref{gnt}), and (\ref{genestsnm}) leads
to%
\begin{equation}
\left\vert t^{m+\ell}g_{n}^{\left(  m\right)  }\left(  t\right)  \right\vert
\leq\frac{\left(  n+\sqrt{n}\right)  ^{m+\ell}}{\sqrt{n+1}}\left(  x+1\right)
^{m+\ell}\exp\left(  \sqrt{n}\left(  1-x\right)  \right)  \quad\text{with\quad
}x=\frac{t}{n+\sqrt{n}}\in\left[  1,\infty\right[  . \label{esttlarge1}%
\end{equation}
The right-hand side in (\ref{esttlarge1}) is maximal for $x=\frac{m+\ell
}{\sqrt{n}}-1$ if $m+\ell\geq2\sqrt{n}$ and for $x=1$ otherwise.

\textbf{Case 1) }$m+\ell\geq2\sqrt{n}$.

For $m+\ell\geq2\sqrt{n}$, the right-hand side in (\ref{esttlarge1}) is
maximal for $x=\frac{m+\ell}{\sqrt{n}}-1$ and we get%
\[
\left\vert t^{m+\ell}g_{n}^{\left(  m\right)  }\left(  t\right)  \right\vert
\leq n^{\frac{m+\ell}{2}}\frac{2^{m+\ell}}{\sqrt{n+1}}\left(  m+\ell\right)
^{m+\ell}\exp\left(  2\sqrt{n}-m-\ell\right)  \leq\left(  2\operatorname*{e}%
\right)  ^{m+\ell}\frac{\left(  m+1\right)  !}{\sqrt{m+1}}n^{\frac{m+\ell
-1}{2}}.
\]

\textbf{Case 2) }$\sqrt{n}\leq m+\ell\leq2\sqrt{n}$.

In this case, the right-hand side in (\ref{esttlarge1}) is maximal for $x=1$
and we arrive, by Stirling's formula and by using $\sqrt{n}\leq m+\ell$, at
the estimate%
\begin{align*}
\left\vert t^{m+\ell}g_{n}^{\left(  m\right)  }\left(  t\right)  \right\vert
&  \leq2^{m+\ell}\frac{\left(  n+\sqrt{n}\right)  ^{m+\ell}}{\sqrt{n+1}}%
\leq\frac{4^{m+\ell}n^{m/2+\ell}}{\sqrt{n+1}}n^{m/2}\\
&  \leq\frac{4^{m+\ell}n^{m/2+\ell}}{\sqrt{n+1}}\left(  m+\ell\right)
^{m}\leq\frac{\left(  4\operatorname*{e}\right)  ^{m+\ell}}{\left(
m+1\right)  ^{1/2}}m!n^{\frac{m-1}{2}+\ell}.
\end{align*}

\textbf{Case 3) }$m+\ell\leq\sqrt{n}$.

Note that%
\[
\left\vert s_{n,m}\left(  n+\delta\right)  \right\vert =\left\vert \tilde
{s}_{n,m}\left(  -\delta\right)  \right\vert =\left\vert \sum_{\ell=0}%
^{m}n^{\ell}p_{\ell,m}\left(  -\delta\right)  \right\vert
\]
with $p_{\ell,m}$ as in (\ref{defplm}). Since $\delta\geq\sqrt{n}\geq m+\ell$,
we get (cf. (\ref{estplmabs})%
\[
\left\vert p_{k,m}\left(  -\delta\right)  \right\vert \leq\dfrac{\left(
m+1\right)  !}{\left(  m-2k\right)  !}\delta^{m-2k}\qquad\forall0\leq k\leq
m/2
\]
and, in turn, (cf. (\ref{defstilde}))%
\begin{align*}
n^{k}\left\vert p_{k,m}\left(  -\delta\right)  \right\vert  &  \leq
\dfrac{\left(  m+1\right)  !}{\left(  m-2k\right)  !}\delta^{m-2k}%
n^{k}=\left(  m+1\right)  !n^{m/2}\frac{\left(  \frac{\delta}{\sqrt{n}%
}\right)  ^{m-2k}}{\left(  m-2k\right)  !}\\
&  \overset{\text{(\ref{repgknm3})}}{\leq}\left\{
\begin{array}
[c]{ll}%
\left(  m+1\right)  \delta^{m} & \delta\geq m\sqrt{n},\\
\left(  m+1\right)  !n^{m/2}\operatorname*{e}^{\delta/\sqrt{n}} & \delta
\geq\sqrt{n}.
\end{array}
\right.
\end{align*}
This implies for $s_{n,m}$ the estimate%
\[
\left\vert s_{n,m}\left(  n+\delta\right)  \right\vert =\left\vert \tilde
{s}_{n,m}\left(  -\delta\right)  \right\vert \leq\left\{
\begin{array}
[c]{ll}%
\left(  m+1\right)  ^{2}\delta^{m} & \delta\geq m\sqrt{n},\\
\left(  m+2\right)  !n^{m/2}\operatorname*{e}^{\delta/\sqrt{n}} & \delta
\geq\sqrt{n}.
\end{array}
\right.
\]
The combination with (\ref{repgknm1}) yields with $t=n+\delta$%
\begin{equation}
\left\vert t^{m+\ell}g_{n}^{\left(  m\right)  }\left(  t\right)  \right\vert
\leq\left\{
\begin{array}
[c]{ll}%
\dfrac{\left(  n+\delta\right)  ^{n+\ell}\operatorname*{e}^{-n-\delta}}%
{n!}\left(  m+1\right)  ^{2}\delta^{m} & \delta\geq m\sqrt{n},\\
\dfrac{\left(  n+\delta\right)  ^{n+\ell}\operatorname*{e}^{-n-\delta}}%
{n!}\left(  m+2\right)  !n^{m/2}\operatorname*{e}^{\delta/\sqrt{n}} &
\delta\geq\sqrt{n}.
\end{array}
\right.  \label{testfinal}%
\end{equation}

\textbf{Case 3a) }$\sqrt{n}\leq\delta\leq m\left(  \sqrt{n}+2\right)  $.

We recall that in this case $m\leq\sqrt{n}$ so that the (generous) estimate
$\delta\leq n+2\sqrt{n}\leq3n$ can be applied.

We employ the second case in the right-hand side of (\ref{testfinal}). First,
we estimate one factor $\left(  n+\delta\right)  ^{\ell}$ by $4n^{\ell}$ and
then observe that the maximum of the remaining expression is taken at
$\delta=\frac{n}{\sqrt{n}-1}$. Thus, we have%
\begin{align*}
\left\vert t^{m+\ell}g_{n}^{\left(  m\right)  }\left(  t\right)  \right\vert
&  \leq4\dfrac{\left(  \frac{n^{\frac{3}{2}}}{\sqrt{n}-1}\right)  ^{n}%
\exp\left(  -\left(  n+\sqrt{n}\right)  \right)  }{n!}n^{m/2+\ell}\left(
m+2\right)  !\\
&  \overset{\text{Stirling}}{\leq}\frac{4}{\sqrt{n+1}}\left(  \frac{\sqrt{n}%
}{\sqrt{n}-1}\right)  ^{n}\exp\left(  -\sqrt{n}\right)  n^{m/2+\ell}\left(
m+2\right)  !.
\end{align*}
Since the function $\left(  \frac{\sqrt{n}}{\sqrt{n}-1}\right)  ^{n}%
\exp\left(  -\sqrt{n}\right)  $ is monotonously decreasing, we get, for
$n\geq2$, the estimate%
\[
\left\vert t^{m+\ell}g_{n}^{\left(  m\right)  }\left(  t\right)  \right\vert
\leq12n^{\frac{m-1}{2}+\ell}\left(  m+2\right)  !.
\]

\textbf{Case 3b) }$m\left(  \sqrt{n}+2\right)  \leq\delta$. We employ the
first case in the right-hand side of (\ref{testfinal}). Since $m\geq1$, the
right-hand side in (\ref{testfinal}) (first case) is maximal for $\delta
_{0}\left(  m\right)  =\frac{m+\ell}{2}\left(  1+\sqrt{1+\frac{4nm}{\left(
m+\ell\right)  ^{2}}}\right)  $. Note that the quantity $\frac{\delta
_{0}\left(  m\right)  }{m}$ is monotonously decreasing with respect to $m$ so
that, for $\ell=0,1$, it holds
\[
\delta_{0}\left(  m\right)  \leq m\left(  1+\sqrt{1+n}\right)  \leq m\left(
\sqrt{n}+2\right)  .
\]
Hence, the right-hand side in (\ref{testfinal}) (first case) has its maximum
at $\delta_{\star}=m\left(  \sqrt{n}+2\right)  $ which is given by (recall
$m\left(  \sqrt{n}+2\right)  \leq3n$)%
\begin{align}
\operatorname*{rhs}  &  :=\dfrac{\left(  n+\delta_{\star}\right)  ^{n+\ell
}\operatorname*{e}^{-n-\delta_{\star}}}{n!}\left(  m+1\right)  ^{2}%
\delta_{\star}^{m}\label{rhsm1}\\
&  \overset{\text{Stirling}}{\leq}\left(  4n\right)  ^{\ell}\left(
\frac{1+\frac{\delta_{\star}}{n}}{\operatorname*{e}^{\delta_{\star}/n}%
}\right)  ^{n}\left(  m+1\right)  ^{2}\delta_{\star}^{m}. \label{rhsm0}%
\end{align}

\textbf{Case 3b}$_{1}$\textbf{) General range: }$m\in\left[  1,\sqrt
{n}\right]  $. We consider the numerator in the right-hand side of
(\ref{rhsm1}) as a function of a free positive variable $\delta_{\star}$ with
maximum at $\delta_{\star}=\ell$. This leads to%
\[
\operatorname*{rhs}\leq\left(  n+1\right)  ^{\ell-1/2}\left(  m+1\right)
^{2}m^{m}3^{m}n^{m/2}\overset{\text{Stirling}}{\leq}\sqrt{2}\left(
3\operatorname*{e}\right)  ^{m}\left(  m+2\right)  !n^{\frac{m-1}{2}+\ell}.
\]

\textbf{Case 3b}$_{2}$\textbf{) Restricted range: }$m\in\left[  2\log
n,\sqrt{n}\right]  $

We will estimate (\ref{rhsm0}) from above. Note that the function $\left(
1+x\right)  \operatorname*{e}^{-x}$ is monotonously decreasing for $x\geq0$.
Since $\frac{\delta_{\star}}{n}\geq\frac{m}{\sqrt{n}}$we get%
\begin{equation}
\operatorname*{rhs}\leq\left(  4n\right)  ^{\ell}\left(  \frac{1+\frac
{m}{\sqrt{n}}}{\operatorname*{e}^{m/\sqrt{n}}}\right)  ^{n}\left(  m+1\right)
^{2}\delta_{\star}^{m}. \label{rhs5}%
\end{equation}

\textbf{Case 3b}$_{2\operatorname*{I}}$\textbf{) }$m\leq\frac{3\sqrt{n}}{4}$.

A Taylor argument for the logarithm implies%
\[
1+x\leq\exp\left(  x-\frac{1}{4}x^{2}\right)  \quad\forall0\leq x\leq3/4
\]
and, in turn,%
\[
\operatorname*{rhs}\leq\frac{\left(  3\operatorname*{e}\right)  ^{m}\left(
m+2\right)  !}{\sqrt{m+1}}\left(  4n\right)  ^{\ell}\left(  n\operatorname*{e}%
\nolimits^{-m/2}\right)  ^{m/2}.
\]
For $m\geq2\log n$ the last bracket is bounded by $1$ and we have proved%
\[
\operatorname*{rhs}\leq\left(  3\operatorname*{e}\right)  ^{m}\left(
m+2\right)  !\left(  4n\right)  ^{\ell}.
\]

\textbf{Case b}$_{2\operatorname*{II}}$\textbf{) }$\frac{3\sqrt{n}}{4}\leq
m\leq\sqrt{n}$.

Note that the bracket in the right-hand side of (\ref{rhs5}) is monotonously
decreasing as a function of $m/\sqrt{n}$ and hence, by choosing $m=3\sqrt
{n}/4$ and using $\delta_{\star}\leq3n$, we end up with%
\[
\operatorname*{rhs}\leq3^{m}\left(  4n\right)  ^{\ell}\left(  m+1\right)
^{2}\left(  c^{n}n^{m}\right)  \quad\text{with\quad}c=\frac{9}{10}.
\]
The last bracket is maximal for $n=\frac{m}{\log\frac{1}{c}}$ so that%
\[
\operatorname*{rhs}\leq\left(  \frac{3}{\log\frac{1}{c}}\right)  ^{m}\left(
4n\right)  ^{\ell}\left(  m+1\right)  ^{2}\operatorname*{e}\nolimits^{-m}%
m^{m}\leq\left(  \frac{3}{\log\frac{1}{c}}\right)  ^{m}\left(  m+2\right)
!\left(  4n\right)  ^{\ell}.
\]
%

%TCIMACRO{\TeXButton{End Proof}{\endproof}}%
%BeginExpansion
\endproof
%EndExpansion

\appendix

\section{Some Auxiliary Estimates}

In this section we first provide an estimate for the approximation of the
function $f\left(  t,x\right)  :=\frac{1}{\sqrt{t^{2}-x^{2}}}$ by its Taylor
polynomial $t^{-1}T_{\mu}\left(  x/t\right)  $ with respect to $x$ around
$x_{0}=0$, where $T_{\mu}$ is as in (\ref{defTmue}).

\begin{proposition}
\label{PropTaylor}Let $0\leq x<t$. Then%
\begin{equation}
\left\vert f\left(  t,x\right)  -t^{-1}T_{\mu}\left(  x/t\right)  \right\vert
\leq\left(  \frac{x}{t}\right)  ^{2\mu}\frac{1}{\sqrt{t^{2}-x^{2}}}.
\label{ftxstatement}%
\end{equation}

\end{proposition}%

%TCIMACRO{\TeXButton{Proof}{\proof}}%
%BeginExpansion
\proof
%EndExpansion
We will prove%
\[
\left\vert 1-\sqrt{1-x^{2}}T_{\mu}\left(  x\right)  \right\vert \leq x^{2\mu
}\qquad\forall0\leq x<1
\]
from which (\ref{ftxstatement}) follows for $t=1$ and, for general $t>0$, by a
simple scaling argument. By using \cite[(5.24.31)]{Hansen75} we obtain%
\[
\frac{1}{\sqrt{1-x^{2}}}=\sum_{k=0}^{\infty}\binom{2k}{k}\left(  \frac{x}%
{2}\right)  ^{2k}\qquad\left\vert x\right\vert <1.
\]
Using $\left(  \sqrt{1-x^{2}}\right)  ^{\prime}=-\frac{x}{\sqrt{1-x^{2}}}$ we
conclude that%
\[
\sqrt{1-x^{2}}=-\sum_{k=0}^{\infty}\frac{1}{\left(  2k-1\right)  }\binom
{2k}{k}\left(  \frac{x}{2}\right)  ^{2k}\qquad\left\vert x\right\vert <1.
\]
Hence,%
\[
1-\sqrt{1-x^{2}}T_{\mu}\left(  x\right)  =1+\left(  \sum_{k=0}^{\infty}%
\frac{1}{\left(  2k-1\right)  }\binom{2k}{k}\left(  \frac{x}{2}\right)
^{2k}\right)  \left(  \sum_{\ell=0}^{\mu-1}\binom{2\ell}{\ell}\left(  \frac
{x}{2}\right)  ^{2\ell}\right)  =1+\sum_{r=0}^{\infty}c_{r,\mu}\left(
\frac{x}{2}\right)  ^{2r}%
\]
with $c_{r,\mu}:=\sum_{\ell=0}^{\min\left\{  r,\mu-1\right\}  }\frac
{1}{2\left(  r-\ell\right)  -1}\binom{2\left(  r-\ell\right)  }{r-\ell}%
\binom{2\ell}{\ell}$. For all $\mu\in\mathbb{N}_{\geq1}$ and $r\in
\mathbb{N}_{0}$ it follows by induction%
\[
c_{r,\mu}=\left\{
\begin{array}
[c]{ll}%
-1 & r=0,\\
0 & 1\leq r<\mu,\\
\frac{\mu}{r}\binom{2\mu}{\mu}\binom{2\left(  r-\mu\right)  }{r-\mu} & \mu\leq
r
\end{array}
\right.
\]
so that%
\begin{align*}
1-\sqrt{1-x^{2}}T_{\mu}\left(  x\right)   &  =\mu\binom{2\mu}{\mu}\left(
\frac{x}{2}\right)  ^{2\mu}\sum_{r=0}^{\infty}\frac{1}{r+\mu}\binom{2r}%
{r}\left(  \frac{x}{2}\right)  ^{2r}\\
&  \leq\mu\binom{2\mu}{\mu}\left(  \frac{x}{2}\right)  ^{2\mu}\sum
_{r=0}^{\infty}\frac{1}{r+\mu}\binom{2r}{r}\left(  \frac{1}{2}\right)  ^{2r}\\
&  =x^{2\mu}.
\end{align*}%
%TCIMACRO{\TeXButton{End Proof}{\endproof}}%
%BeginExpansion
\endproof
%EndExpansion

\begin{proposition}
\label{Propgnderv}For any $m,n\in\mathbb{N}_{0}$ and $0\leq2r\leq n-1$, it
holds%
\[
\int_{x}^{\infty}t^{m-1-2r}g_{n}^{\left(  m\right)  }\left(  t\right)
dt=-\frac{\left(  n-1-2r\right)  !}{n!}\left(  \frac{d}{dx}\right)
^{2r}\left(  x^{m}g_{n-1-2r}^{\left(  m-1-2r\right)  }\left(  x\right)
\right)  \quad\forall x>0\text{.}%
\]

\end{proposition}%

%TCIMACRO{\TeXButton{Proof}{\proof}}%
%BeginExpansion
\proof
%EndExpansion
Since the limit $x\rightarrow\infty$ of both sides converges to zero due to
the exponential decay of $g_{n}^{\left(  k\right)  }\left(  x\right)  $ it
sufficient to prove that the derivatives of both sides coincide. We set%
\begin{align*}
L_{r}^{m,n}\left(  x\right)   &  :=-x^{m-1-2r}g_{n}^{\left(  m\right)
}\left(  x\right) \\
R_{r}^{m,n}\left(  x\right)   &  :=-\frac{\left(  n-1-2r\right)  !}{n!}\left(
\frac{d}{dx}\right)  ^{2r+1}\left(  x^{m}g_{n-1-2r}^{\left(  m-1-2r\right)
}\left(  x\right)  \right)
\end{align*}
and prove $L_{r}^{m,n}\left(  x\right)  =R_{r}^{m,n}\left(  x\right)  $ by
induction over $r$.

\textbf{Start of induction: }$r=0$. Then,%
\[
L_{0}^{m,n}\left(  x\right)  =-x^{m-1}g_{n}^{\left(  m\right)  }\left(
x\right)  .
\]
By $m-$times differentiating the relation $g_{n}\left(  x\right)  =\frac{x}%
{n}g_{n-1}\left(  x\right)  $ we obtain%
\begin{align*}
L_{0}^{m,n}\left(  x\right)   &  =-x^{m-1}g_{n}^{\left(  m\right)  }\left(
x\right)  =\frac{-x^{m-1}}{n}\left(  xg_{n-1}^{\left(  m\right)  }\left(
x\right)  +mg_{n-1}^{\left(  m-1\right)  }\left(  x\right)  \right) \\
&  =-\frac{1}{n}\frac{d}{dx}\left(  x^{m}g_{n-1}^{\left(  m-1\right)  }\left(
x\right)  \right)  =R_{0}^{m,n}\left(  x\right)  .
\end{align*}

\textbf{Induction step: }We assume that $L_{k}^{m,n}=R_{k}^{m,n}$ for $0\leq
k\leq r$. Then%
\begin{align*}
R_{r+1}^{m,n}\left(  x\right)   &  :=-\frac{\left(  n-3-2r\right)  !}%
{n!}\left(  \frac{d}{dx}\right)  ^{2r+1}\left(  \frac{d}{dx}\right)
^{2}\left(  x^{m}g_{n-3-2r}^{\left(  m-3-2r\right)  }\left(  x\right)  \right)
\\
&  =\frac{m\left(  m-1\right)  R_{r}^{m-2,n-2}\left(  x\right)  +2mR_{r}%
^{m-1,n-2}\left(  x\right)  +R_{r}^{m,n-2}\left(  x\right)  }{n\left(
n-1\right)  }.
\end{align*}
By using the induction assumption and several times the recurrence relation%
\[
g_{n}\left(  x\right)  =\frac{x}{n}g_{n-1}\left(  x\right)  \quad\text{which
implies\quad}g_{n}^{\left(  k\right)  }\left(  x\right)  =\frac{x}{n}%
g_{n-1}^{\left(  k\right)  }\left(  x\right)  +\frac{k}{n}g_{n-1}^{\left(
k-1\right)  }\left(  x\right)
\]
we obtain%
\begin{align*}
R_{r+1}^{m,n}\left(  x\right)   &  =\frac{m\left(  m-1\right)  L_{r}%
^{m-2,n-2}\left(  x\right)  +2mL_{r}^{m-1,n-2}\left(  x\right)  +L_{r}%
^{m,n-2}\left(  x\right)  }{n\left(  n-1\right)  }\\
&  =-\frac{m\left(  m-1\right)  x^{m-3-2r}g_{n-2}^{\left(  m-2\right)
}\left(  x\right)  +2mx^{m-2-2r}g_{n-2}^{\left(  m-1\right)  }\left(
x\right)  +x^{m-1-2r}g_{n-2}^{\left(  m\right)  }\left(  x\right)  }{n\left(
n-1\right)  }\\
&  =-\frac{mx^{m-3-2r}}{n}g_{n-1}^{\left(  m-1\right)  }-\frac{x^{m-2-2r}}%
{n}g_{n-1}^{\left(  m\right)  }\\
&  =-x^{m-3-2r}g_{n}^{\left(  m\right)  }=L_{r+1}^{m,n}\left(  x\right)  .
\end{align*}%
%TCIMACRO{\TeXButton{End Proof}{\endproof}}%
%BeginExpansion
\endproof
%EndExpansion

\begin{lemma}
\label{LemDoubleFactorial}For $m\geq0$ and $0\leq k\leq m$, it holds%
\[
\left(  \frac{3}{5}\right)  ^{k}\sqrt{\frac{m!}{\left(  m-k\right)  !}}%
\leq\frac{m!!}{\left(  m-k\right)  !!}\leq2^{k}\sqrt{\frac{m!}{\left(
m-k\right)  !}}%
\]

\end{lemma}%

%TCIMACRO{\TeXButton{Proof}{\proof}}%
%BeginExpansion
\proof
%EndExpansion
The case $m=0$ and the case $k=0$ are trivial. Let $k=1$. The formula is easy
to check for $m=1,2$ and we restrict in the following to $m\geq3$.

We use Stirling's formula in the form%
\[
n!=C_{n}n^{n+\frac{1}{2}}\exp\left(  -n\right)  \quad\text{with\quad}%
C_{n}=\sqrt{2\pi}\exp\left(  \theta/\left(  12n\right)  \right)  \text{ for
}n\in\mathbb{N}_{\geq1}\text{ and }\theta\in\left]  0,1\right[
\]
so that $\frac{C_{n}}{\sqrt{2\pi}}\in\left]  1,\operatorname*{e}%
^{1/12}\right[  $.

For $m=2r+1$, we get with $r\geq1$%
\begin{align*}
\frac{\left(  2r+1\right)  !!}{\left(  2r\right)  !!}  &  =\frac{\left(
2r+1\right)  !}{4^{r}\left(  r!\right)  ^{2}}=\frac{C_{m}}{C_{r}^{2}}%
\frac{\left(  2r+1\right)  ^{2r+\frac{3}{2}}\exp\left(  -2r-1\right)  }%
{4^{r}r^{2r+1}\exp\left(  -2r\right)  }\\
&  =\frac{C_{m}}{C_{r}^{2}\operatorname*{e}}2\left(  1+\frac{1}{2r}\right)
^{2r+1}\sqrt{2r+1}\left\{
\begin{array}
[c]{ll}%
\leq\frac{\operatorname*{e}^{1/12-1}}{\sqrt{2\pi}}2\left(  \frac{3}{2}\right)
^{3}\sqrt{m} & \leq\frac{11}{10}\sqrt{m},\\
\geq\frac{2}{\sqrt{2\pi}\operatorname*{e}^{1/6}}\sqrt{m} & \geq\frac{3}%
{5}\sqrt{m}.
\end{array}
\right.
\end{align*}
For $m=2r$ and $r\geq2$, we get%
\begin{align*}
\frac{\left(  2r\right)  !!}{\left(  2r-1\right)  !!}  &  =\frac{\left(
2r\right)  !!\left(  2r-2\right)  !!}{\left(  2r-1\right)  !}=\frac
{2^{2r-1}r!\left(  r-1\right)  !}{\left(  2r-1\right)  !}=\frac{C_{r}%
C_{r-1}2^{2r-1}r^{r+1/2}\left(  r-1\right)  ^{r-1/2}\operatorname*{e}^{1-2r}%
}{C_{m-1}\left(  2r-1\right)  ^{2r-1/2}\operatorname*{e}^{1-2r}}\\
&  =\sqrt{m-1}\frac{C_{r}C_{r-1}}{2C_{m-1}}\left(  \frac{r}{r-1/2}\right)
^{r+1/2}\left(  \frac{r-1}{r-1/2}\right)  ^{r-1/2}\\
&  \left\{
\begin{array}
[c]{c}%
\leq\sqrt{\frac{\pi\left(  m-1\right)  }{2}}\operatorname*{e}^{-1/3}\left(
\frac{4}{3}\right)  ^{5/2}\leq2\sqrt{m-1}\leq2\sqrt{m},\\
\geq\sqrt{\frac{\pi\left(  m-1\right)  }{2}}\operatorname*{e}^{\frac{5}{12}%
}\left(  \frac{2}{3}\right)  ^{3/2}\geq\frac{11}{10}\sqrt{m-1}\geq\frac{4}%
{5}\sqrt{m}.
\end{array}
\right.
\end{align*}
In summary, we have proved%
\[
\frac{3}{5}\sqrt{m}\leq\frac{m!!}{\left(  m-1\right)  !!}\leq2\sqrt{m}.
\]
From this the assertion follows by induction.%
%TCIMACRO{\TeXButton{End Proof}{\endproof}}%
%BeginExpansion
\endproof
%EndExpansion

\bibliographystyle{abbrv}
\bibliography{nlailu}

\end{document}